\documentclass[a4paper,12pt]{article}

\usepackage[latin1]{inputenc}
\usepackage[francais,english]{babel}
\usepackage{graphicx,url}
\usepackage{latexsym,amssymb,amsmath,amsthm}
\usepackage[T1]{fontenc}

\theoremstyle{definition}

\theoremstyle{remark}

\numberwithin{equation}{section}

\begin{document}

\title{Alberti's letter counts}
\author{Bernard Ycart\\
Laboratoire Jean Kuntzmann\\ 
Universit\'e Joseph Fourier and CNRS UMR 5224\\ 
51 rue des Math\'ematiques 38041 Grenoble cedex 9, France \\
\texttt{Bernard.Ycart@imag.fr}}
\maketitle
\begin{abstract}
Four centuries before modern statistical linguistics was born, Leon
Battista Alberti (1404--1472) compared the
frequency of vowels in Latin poems and orations, making the first
quantified observation of a stylistic difference ever. Using a corpus of 20
Latin texts (over 5 million letters), Alberti's observations are
statistically assessed. Letter counts prove that poets used
significantly more \emph{a}'s, \emph{e}'s, and \emph{y}'s, whereas orators
used more of the other vowels. The sample sizes needed to
justify the assertions are studied, and proved to be within
reach for Alberti's scholarship.
\end{abstract}

{\small

}
\section{Introduction}
\label{intro}
August 18, 1851 is the date of birth generally agreed upon for
statistical linguistics \cite{Williams1956,Bailey1969}. 
On that day, Augustus de Morgan (1806--1871) sent a letter to a
friend, suggesting
statistical counts to settle 
authorship disputes \cite{Morgan1851}. This later inspired the
first statistical study of the kind by T. C. Mendenhall (1841--1924)
 \cite{Mendenhall1887}. 
\vskip 2mm\noindent
Four centuries before, Leon Battista
 Alberti (1404--1472), had made a curious remark on the frequency of vowels
 among poets and orators. The first observation of the kind ever, it
 has remained unnoticed since. That remark, and its statistical
 verification, are the subject of this note.
\vskip 2mm\noindent
Alberti, the well known polymath of Italian
Renaissance, has been hailed as a pioneer in quite different domains
\cite{Grafton2000}:
architecture, grammar, mathematics, etc. His
``De componendis cifris'', written in 1466
or 1467, stands out as the
first western text in the history of cryptology \cite{Kahn1973}. The
latin text is available online \cite{Alberti1466}; several
translations have been published, among which we shall mainly quote
\cite{Alberti2010}; we have also used an early Italian translation by
C. Bartoli \cite{Bartoli1568}, and M. Furno's French translation
\cite{Alberti2000}. Our focus is on the following passage of 
section IV.
\begin{quotation}
Sic enim adnotasse videor apud poetas
vocales a consonantibus numero superari
non amplius quam ex octava; apud
rhetores vero non excedere consonantes
ferme ex proportione quam sesqui\-tertiam
nuncupant. Nam si fuerint quidem
connumeratae in unumque collectae
omnes istius generis paginae vocales
numero puta tricentarum, reliquarum
omnium consonantium numerus una
coadiunctus erit fere quadringentarum.
\end{quotation}
Here is K. Williams' translation \cite[p.~173]{Alberti2010}.
\begin{quotation}
From my calculations, it turns out that in the case of poetry, the
number of consonants exceeds the number of vowels by no more than an
octave, while in the case of prose the consonants do not usually
exceed the vowels by a ratio greater than a sesquialtera. If in
fact we add up all the vowels on a page, let's say there are three
hundred, the overall sum of the consonants will be four hundred.
\end{quotation}
The interpretation of Alberti's assertions in this passage,  has been
debated. Contrarily to note 3 in \cite[p.~173]{Alberti2010}, we
believe that the correct translation of ``ex octava'' in the first
sentence is indeed
``one eighth'', as in several other translations. For convenience, the
proportion of vowels for poets will be denoted by $P$, being aware that
what precisely is meant by
``poets'' and ``vowels'', is far from clear at this point. We believe
that Alberti's assertion can be mathematically 
translated into $(1-P)-P<1/8$, or else
$P>7/16$. Alberti opposes poets to ``rhetores'', i.e. orators,
for which the proportion of vowels against consonants 
is said to be ``sesquitertia'',
translated by  \cite{Bartoli1568} as ``del terzo piu'', or else
four against three. If we denote by $R$ the proportion of vowels for orators
(with the same precautions as before), what is stated is  
$(1-R)/R < 4/3$, or else $R>3/7$. The difference between $7/16=
43.75\%$ and $3/7=42.86\%$ is smaller than $1\%$: how did Alberti
observe and justify such a small difference? The last sentence of the
passage contains a hint. At first glance, saying that for 300 vowels,
400 consonants are counted, seems redundant, once the proportion has
been set to 3:4. Yet ``De componendis cifris'' is a
rather short and concentrated text, and no sentence is
superfluous. Mentioning a page of 700 letters might have been
Alberti's way of indicating that his observations were supported by
counts on large enough sets of letters, 
i.e. the modern notion of sample size. 
\vskip 2mm\noindent
Admittedly, Alberti is far from being the first to have counted
letters. Six centuries before, Al Kind\={\i} (ca. 801--873)
had written a
``Manuscript on deciphering cryptographic messages''
\cite{AlKadi1992,Mrayatietal2002}. Remarkably similar in its
organization and 
contents to Alberti's ``De componendis cifris'',  Al Kind\={\i}'s
treatise contains a table of letter frequencies, much more detailed
than Alberti's observations. It also contains quite accurate
linguistic observations on poetry and prose. Similar counts and
remarks can also be found in later Arab treatises, in particular that
of ibn Dunayn\={\i}r (1187--1229) \cite{AlKadi1992,Mrayatietal2005}. 
Nevertheless, all Arab treatises, as well as all western
writings after Alberti until the middle of the
19\textsuperscript{th} century, viewed letter frequencies as
a \emph{characteristic of a language}. Alberti was the first to assert
that a stylistic difference (poetry vs. oration) in the same language,  
could yield a quantifiable difference in letter frequencies.
\vskip 2mm\noindent
Our first objective was to investigate whether Alberti's
observations could be justified, by modern statistical standards. 
We gathered a corpus of classical latin texts, from 10 different
poets and as many different orators, totalling
over 5 million letters. The letter frequency analysis encountered
the difficulty of distinguishing the uses of \emph{u} and \emph{i},
as vowel or  consonant. A counting algorithm was proposed, based on
classical Latin grammar.  Using that grammatical algorithm, the vowel counts
did show a clear difference between poets and orators: poets
use consistently more \emph{a}'s, \emph{e}'s and \emph{y}'s, whereas 
orators use more of the other vowels. However, the total 
percentages of vowels among 
poets ($42.92\%$) and among orators ($43.14\%$) were not significantly
different, the percentage among orators being even slightly
higher. Another counting algorithm did meet
Alberti's observations: it consisted in
counting as consonants all \emph{i}'s before another vowel. 
Even though it cannot be asserted that our counting algorithm matches
his, we believe that our
results support Alberti's observations, and that he had indeed
detected the quantitative difference between poets and orators.
\vskip 2mm\noindent 
Our second objective was to calculate how many pages of
700 letters would have been necessary to statistically assess
Alberti's observations. We translated them into 3 statistical tests,
and computed the number of pages necessary to raise the power of the
tests to $95 \%$, at threshold $5 \%$. Detailed results will be
reported for the 3 tests, using both the classical Bernoulli model,
and a Monte-Carlo study on the corpus. One of our conclusions 
is that, randomly selecting 2 pages in each of the
20 texts of our corpus would suffice to ensure a power of $95 \%$ for
the test of $P=R$ against $P>R$.
\vskip 2mm\noindent
Alberti is known as a
sholar of immense culture, a lover of 
litterature, grammar and alphabets: see \cite{Patota1999} and in particular 
note 127 p.~\textsc{XLIII} of \cite{Alberti2003}.  
He might well have counted letters in 40 pages of poems and
orations: he could definitely have proved his assertions\ldots~ had he
known statistics. 
As he says in section III 
\cite[p.~172]{Alberti2010}:
\begin{quotation}
I went about this with not indifferent industry or care, reflecting
again and again on the elements of writing, and investigating
intensely until I had clear in mind some fundamental concepts, which
now the most brilliant minds will acknowledge as having contributed to
understanding the whole question of ciphers.
\end{quotation}
\vskip 2mm\noindent
Two sections follow; in the first one,
our corpus of texts will be presented, and
different ways of counting vowels will be discussed. The vowel
frequency difference between poets and orators will be
statistically assessed. The sample size problem is treated in section
\ref{stats}. Three tests will be defined, and their powers computed
using both the Bernoulli model, and random pages extracted from our corpus.
The free statistical software R was used for letter counts and
simulation experiments (\cite{R,Gries2009}. The corpus, the R script
of functions, and a user manual have been made available online as a
compressed 
file\footnote{http://ljk.imag.fr/membres/Bernard.Ycart/publis/llc.tgz}.
\section{Corpus and vowel counts}
\label{corpus}
\begin{figure}[!ht]
\centerline{
\includegraphics[width=12cm]{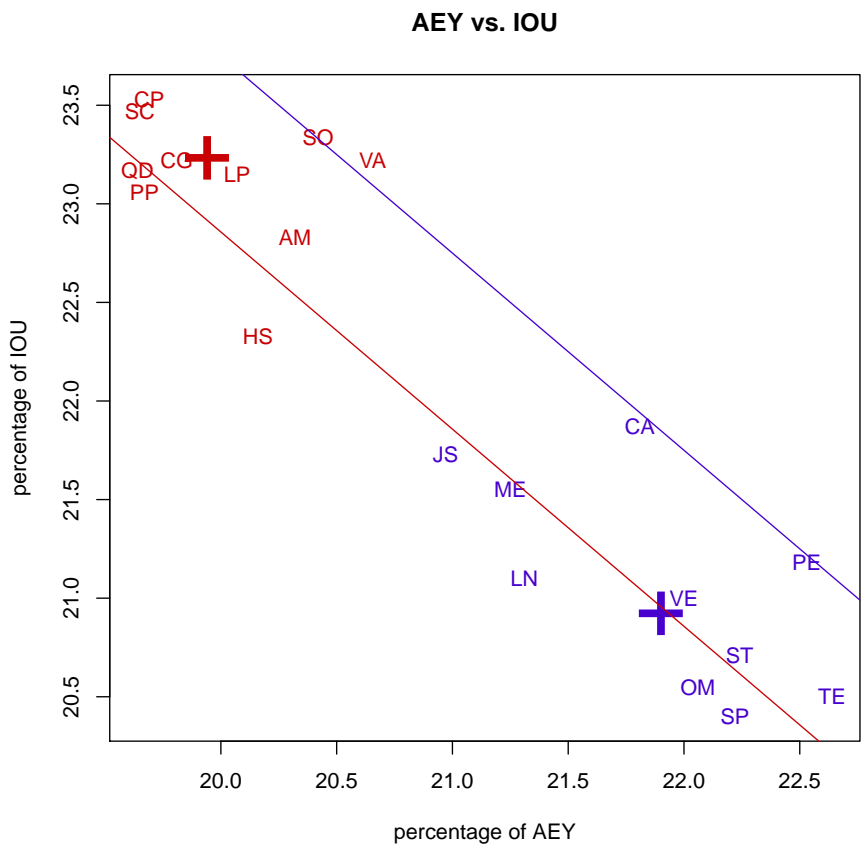}
} 
\caption{Vowel proportions in grammatical counts. Each text is
  represented by its code, color blue for poets, red for orators. The
  abscissa of a text is the cumulated percentage of \emph{a}'s,
  \emph{e}'s, and \emph{y}'s
the ordinate is the cumulated percentage of the other
  vowels. The results on the global corpuses of poets and orators are
  represented by a blue and a red cross. Poets use significantly more
  \emph{a,e,y}, and less of the other vowels. The two lines
  correspond to Alberti's assertion: $x+y=7/16$ (in blue for poets),
$x+y=3/7$ (in red for orators). }
\label{fig:AvI1}
\end{figure}

\begin{figure}[!ht]
\centerline{
\includegraphics[width=12cm]{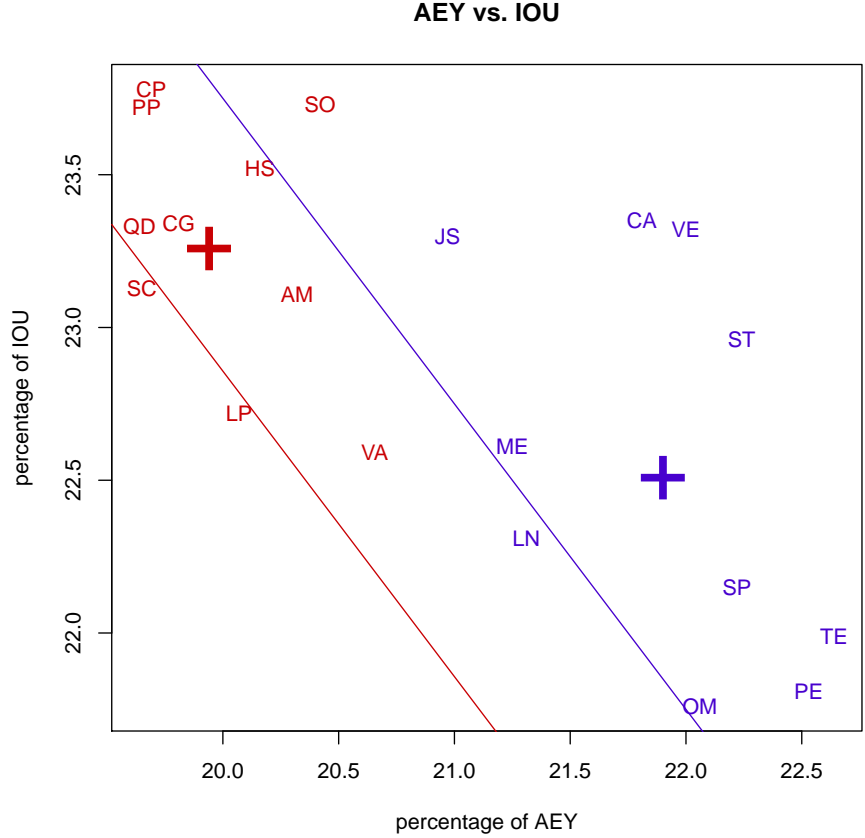}
} 
\caption{Vowel proportions, counting \emph{i} as consonant before
  \emph{a,e,o,u}, with the same conventions as in Figure
  \ref{fig:AvI1}. Alberti's assertions are verified: poets are above the
  blue line on average, orators above the red line.}
\label{fig:AvI2}
\end{figure}

Alberti's encyclopaedic knowledge makes it a hopeless task to guess
from which texts he could have made his letter counts.  
We selected ten of the most famous classic Latin poets 
and orators, and copied their
texts from
``The Latin Library''\footnote{http://www.thelatinlibrary.com}.
Each text was edited to remove non ascii characters, and also Roman
numerals that could have biased letter counts. The texts were
given a two letter code for further reference.
\begin{itemize}
\item Poets (total: 2617488 letters)
\begin{itemize}
\item CA:
Catullus, Poems (71747 l.)\\
Gaius Valerius Catullus, ca. 84 BC -- ca. 54 BC
\item JS:
Juvenal, Saturae (142645 l.)\\
Decimus Junius Juvenalis, 1\textsuperscript{st}-2\textsuperscript{nd} 
century AD
\item LN:
Lucretius, De rerum Natura (274355 l.)\\
Titus Lucretius Carus, ca. 99 BC -- ca. 55 BC
\item ME:
Martial, Epigrams (299099 l.)\\
Marcus Valerius Martialis, 40 AD -- ca. 103 AD
\item OM:
Ovid, Metamorphoses (446848 l.)\\
Publius Ovidius Naso, 43 BC -- ca. 17 AD
\item PE:
Propertius, Elegiae (134719 l.)\\
Sextus Propertius, ca. 50 BC -- ca. 15 BC
\item SP:
Silius Italicus, Punica (449903 l.)\\
Tiberius Catius Asconius Silius Italicus, ca. 28 BC -- ca. 103 BC
\item ST: 
Statius, Thebaid (365326 l.)\\
Publius Papinius Statius, ca. 45 AD -- ca. 96 AD
\item TE:
Tibullus, Elegiae (66062 l.)\\
Albius Tibullus, ca. 55 BC -- 19 BC
\item VE:
Virgil, Aeneid (366784 l.)\\
Publius Vergilius Maro, 70 BC -- 19 BC
\end{itemize}
\item Orators (total: 2570344 letters)
\begin{itemize}
\item AM: 
Apuleius, Metamorphoses (332617 l.)\\
Lucius Apuleius, ca. 125 AD -- ca. 180 AD
\item CG:
Caesar, De Bello Gallico (317056 l.)\\
Gaius Julius Caesar, 100 BC -- 44 BC
\item CP:
Cicero, Catilinarian and Philippics (370521 l.)\\
Marcus Tullius Cicero, 106 BC -- 43 BC
\item HS:
Horace, Sermones   (77859 l.)\\
Quintus Horatius Flaccius, 65 BC -- 8 BC
\item LP:
Lactantius, De Mortibus Persecutorum (65616 l.)\\
Lucius Caecilius Firmianus Lactantius, ca. 240 AD -- ca. 320 AD
\item PP:
Pliny the younger, Panegyricus (112992 l.)\\  
Gaius Plinius Caecilius Secundus, 61 AD -- 112 AD
\item QD:
Quintillian, Declamatio Major (390261 l.)\\
Marcus Fabius Quitilianus,  35 AD -- 100 AD
\item SC:
Seneca the Elder, Controversiae (499280 l.)\\
Marcus Annaeus Seneca, ca. 54 BC -- ca. 39 AD
\item SO:
Sallust, Orations (50676 l.)\\
Gaius Sallustus Crispus, 86 BC -- ca. 35 BC
\item VA:
Vitruvius, De Architectura (353466 l.)\\
Marcus Vitruvius Pollio, ca. 75 BC -- ca. 15 BC
\end{itemize}
\end{itemize}
Our choices of Caesar and Vitruvius in the list of orators are
questionable. Even though ``De Bello Gallico'' is not
explicitly written as a speech, it can be argued that Caesar's style
owed much to his expertise as a politician.
Not being known as an orator, Vitruvius was included
because his ``De Architectura'' was an important source of inspiration for 
Alberti \cite{Grafton2000}. Statistical evidence showed that both texts
behaved similarly to those of other orators regarding letter counts.
\vskip 2mm\noindent
If vowels are counted on the basis of the occurrence of characters
\emph{a}, \emph{e}, \emph{i}, \emph{o}, \emph{u}, \emph{y},
Alberti's assertions can hardly be
justified: the vowel proportions are far from the announced
values ($45.76\%$ for poets and $45.15\%$ for orators instead of 
$43.75\%$ and $42.86\%$). Moreover, the
vowel percentage is not consistently higher for poets. However, it is
a well known particularity of Latin grammar that the letters \emph{i}
and \emph{u} can be used as vowels or as consonants. Alberti was
perfectly aware of that; here are some quotations from sections IV, VI
and VIII of \cite{Alberti2010}.
\begin{quotation}
After the vowel O is sometimes found the I and also the U, but this is
rather rare, while more frequently the O is followed by the U used as
a consonant, as in  `ovem'.

[\ldots]

In fact, neither I or V used as consonants ever follow the vowels in a
syllable, and neither does Q.

[\ldots]

In monosyllabic words the vowels can be followed by all the consonants
except F, G, P, Q, I, and V. 

[\ldots]

In contrast, within a word, there are some consonants that never
appear in the next syllable after V and the consonant I.
\end{quotation}
The usage in rendering vocalic [i], [u] and
consonantal [j], [w] into letters has varied, even for classical Latin.
In our corpus, the letter \emph{i} (lower or upper case) always
denoted either sounds [i] and [j]. Lower case \emph{v}'s were
rather rare, and upper case \emph{U}'s never appeared. 
In order to define what might be called a ``grammatical
count'', we have used both Alberti's own remarks, together with the
rules that can be found in most Latin grammars, such as
\cite{Tafel1860,Blair1874,Allen1978}. Here are the counting rules that
were implemented.
\begin{itemize}
\item \emph{I} or \emph{i} was counted as a consonant:
\begin{itemize}
\item at the beginning of a word before a vowel,
\item after prefixes \emph{ad, ab, conj, ex}, and before another
    vowel,
\item between two vowels.
\end{itemize}
\item \emph{U} or \emph{u} was counted as a consonant:
\begin{itemize}
\item at the beginning of a word before a vowel,
\item after \emph{q} and \emph{g}, and before a vowel,
\item between two vowels,
\end{itemize}
\item \emph{V} or \emph{v} was counted as a vowel before a consonant.
\end{itemize} 
Our choice to count \emph{qu} as two consonants instead of one is
disputable: see \cite[p.~20]{Tafel1860}, \cite[p.~44]{Blair1874},
\cite[p.~16]{Allen1978}. We believe it is supported by the following
remark in section VI, \cite[p.~175]{Alberti2010}:
\begin{quotation}
On the other hand, when it was established that the U is combined with
the Q, it doesn't appear to have been taken into account that the U
itself is implicit in this letter Q, so it sounds like KU.
\end{quotation}
Many exceptions to the above rules
are recorded in grammars; for instance, \emph{i} often has the
consonantal sound [j] after a consonant, before a vowel other than
\emph{i}. No implementation of a
count including exceptions could be attempted. 
Vowel percentages in our grammatical count 
are reported in Table \ref{tab:POt}. The difference between poets and
orators is clear: poets use more \emph{a}'s, \emph{e}'s, and \emph{y}'s and less
of the other vowels. 
In order to assess statistical significance, Student's T-test 
 was applied to paired lines of Table
\ref{tab:POt}: for each letter, the 10 proportions among poets were
compared  to the 10 proportions among orators, and the one-sided
p-value was returned.
\begin{itemize}
\item percentage of \emph{A}: larger for poets, p-value $ = 4.12\times 10^{-7}$,
\item percentage of \emph{E}: larger for poets, p-value $ = 2.10\times 10^{-3}$,
\item percentage of \emph{I}: smaller for poets, p-value $ = 1.98\times 10^{-9}$,
\item percentage of \emph{O}: smaller for poets, p-value $ = 3.01\times 10^{-2}$,
\item percentage of \emph{U}: smaller for poets, p-value $ = 1.95\times 10^{-2}$,
\item percentage of \emph{Y}: larger for poets, p-value $ = 6.51\times 10^{-5}$.
\end{itemize}
\setlength{\tabcolsep}{3pt}
\begin{table}
\begin{center}
\begin{tabular}{|c|rrrrrrrrrr|r|}
\hline
\multicolumn{12}{|c|}{Poets}\\\hline
&CA&JS&LN&ME&OM&PE&SP&ST&TE&VE&All\\\hline
A&     9.38& 9.62& 9.05& 9.60& 9.52&10.71&10.38& 9.83& 9.92& 9.75& 9.77\\
E&    12.15&11.22&12.24&11.44&12.35&11.58&11.64&12.18&12.60&12.06&11.95\\
I&     9.89& 9.80& 9.70& 9.98& 9.48& 9.69& 8.77& 9.05& 9.25& 9.31& 9.39\\
O&     5.31& 5.28& 5.21& 5.47& 5.19& 5.71& 5.15& 5.03& 5.16& 4.99& 5.20\\
U&     6.67& 6.65& 6.19& 6.10& 5.88& 5.78& 6.48& 6.63& 6.09& 6.70& 6.32\\
Y&     0.28& 0.13& 0.02& 0.21& 0.19& 0.24& 0.20& 0.23& 0.12& 0.19& 0.18\\\hline
&43.76&42.78&42.45&42.83&42.69&43.77&42.74&43.13&43.21&43.18&42.92\\\hline\hline
\multicolumn{12}{|c|}{Orators}\\\hline
&AM&CG&CP&HS&LP&PP&QD&SC&SO&TO&All\\\hline
A&     8.75& 8.02& 8.11& 8.68& 8.39& 8.28& 8.05& 8.28& 9.08& 8.29& 8.28\\
E&    11.50&11.79&11.56&11.43&11.64&11.37&11.58&11.27&11.32&11.72&11.51\\
I&    10.85&10.99&11.31&10.20&11.23&11.05&11.43&11.77&11.45&10.95&11.22\\
O&     5.35& 5.72& 5.68& 5.36& 5.05& 5.27& 5.44& 5.59& 5.30& 5.77& 5.52\\
U&     6.63& 6.51& 6.54& 6.77& 6.87& 6.74& 6.30& 6.11& 6.59& 6.96& 6.51\\
Y&     0.07& 0.00& 0.02& 0.05& 0.04& 0.02& 0.01& 0.10& 0.02& 0.02&0.04\\\hline
&43.21&43.17&43.34&42.52&43.30&42.82&42.86&43.18&43.88&43.68&43.14\\\hline
\end{tabular}
\caption{Percentage of vowels in the texts of
  our corpus, counting \emph{i}'s, \emph{u}'s, and \emph{v}'s
  according to grammatical rules.}
\label{tab:POt}
\end{center}
\end{table}
\vskip 2mm\noindent
Figure \ref{fig:AvI1} shows a two-dimensional
representation of the same data, with percentages of \emph{A,E,Y} on the x-axis,
percentages of \emph{I,U,O} on the y-axis. The difference between poets
and orators is quite visible. 
On the same graphic, the two lines $x+y=7/16$ (red) and
$x+y=3/7$ (blue) correspond to Alberti's assertions. Orators are
indeed above the red line on average, but poets are below the blue line. No
clear difference on vowel counts
can be seen, and the global percentage of vowels 
is even slightly smaller for orators than for
poets: the grammatical count does not match Alberti's
observations.  
\vskip 2mm\noindent
The difference between Alberti's and grammatical counts obviously
comes from the letter \emph{i} as a vowel, which remains too
frequent among orators. Another count was proposed, that consisted in
counting as consonant all \emph{i}'s before \emph{a,e,o,u}, counting all
\emph{u}'s as vowels and \emph{v}'s as consonants. 
Table \ref{tab:PO2} shows the total vowel counts, and Figure
\ref{fig:AvI2} the graphical representation, with the same conventions
as Figure \ref{fig:AvI1}. This time, Alberti's observations are
verified: poets are above the blue line, and orators above the red
one, on average. The following comparisons are statistically
significant, using again Student's T-test, and denoting by $P$ and $R$ the
proportions of vowels among poets and orators respectively.
\begin{itemize}
\item $P>7/16$: p-value $ = 1.16\times 10^{-3}$
\item $R>3/7$: p-value $ = 2.24\times 10^{-3}$
\item $P>R$: p-value $ = 1.79\times 10^{-5}$
\end{itemize}
\begin{table}
\begin{center}
\begin{tabular}{|rrrrrrrrrr|r|}
\hline
\multicolumn{11}{|c|}{Poets}\\\hline
CA&JS&LN&ME&OM&PE&SP&ST&TE&VE&All\\\hline
45.15&44.28&43.62&43.86&43.82&44.35&44.37&45.20&44.63&45.32&44.41\\\hline\hline
\multicolumn{11}{|c|}{Orators}\\\hline
AM&CG&CP&HS&LP&PP&QD&SC&SO&TO&All\\\hline
43.43&43.15&43.49&43.69&42.79&43.38&42.97&42.79&44.15&43.24&43.19\\\hline
\end{tabular}
\caption{Percentage of vowels in the texts of
  our corpus, counting \emph{i} as consonant before \emph{a,e,o,u}.}
\label{tab:PO2}
\end{center}
\end{table}
\vskip 2mm\noindent
Of course, verifying Alberti's conclusions with our second counting algorithm,
 does not prove that his counts would have matched ours. Alberti's
 way of counting will always remain a mystery. Yet we believe that our
 results support his observations, and that he had indeed detected the
 quantitative difference between poets and orators.

\section{Sample sizes}
\label{stats}
\begin{figure}[!ht]
\centerline{
\includegraphics[width=14cm]{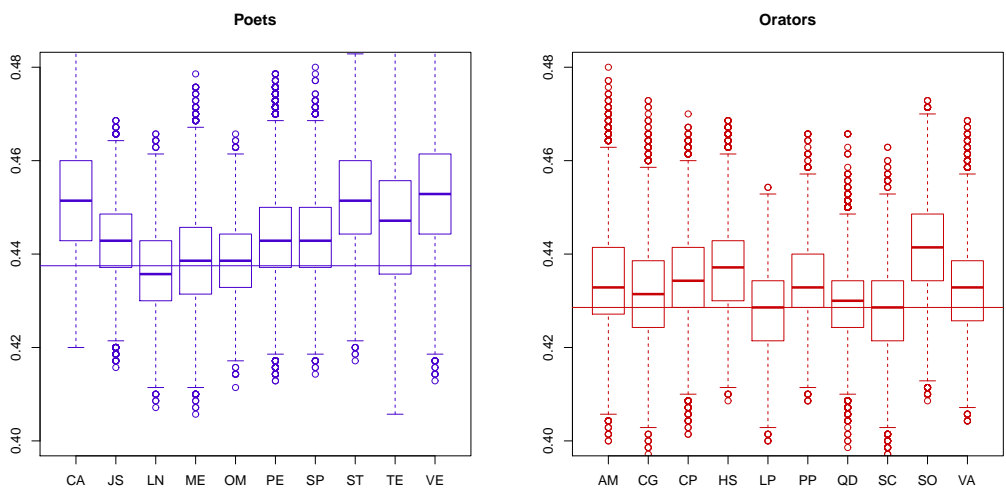}
} 
\caption{Boxplots for the proportions of vowels in 10000 random pages
  of each text. The horizontal lines correspond to Aberti's
  observations: 7/16 for poets, 3/7 for orators.}
\label{fig:PO}
\end{figure}

Our focus in this section is on the last sentence of Alberti's
assertions, about counting vowels in a page of 700 letters: could that
sample size be sufficient to support Alberti's statistical conclusions? If not,
how many pages should be counted? We shall first answer the question
theoretically on the classical Bernoulli model, then show that the minimal
sample sizes are somewhat lower using random pages from the actual corpus.
All statistical techniques
used here are quite classical: see \cite{Gries2009} as a general
reference. 
\vskip 2mm\noindent
We rely upon the observations made on our corpus using the
second counting algorithm, that matches Alberti's observations (Table
\ref{tab:PO2} and Figure \ref{fig:AvI2}). Denoting by $P$ and $R$ the
respective proportions of vowels among poets and among orators,  we
consider the following tests, were $\mathcal{H}_0$ and $\mathcal{H}_1$
denote as usual the null hypothesis and the alternative. 
\begin{enumerate}
\item[[T1\!\!]] $\mathcal{H}_0~: P=7/16$ against $\mathcal{H}_1: P>7/16$,
\item[[T2\!\!]] $\mathcal{H}_0~: R=3/7$ against $\mathcal{H}_1: R>3/7$,
\item[[T3\!\!]] $\mathcal{H}_0~: P=R$ against $\mathcal{H}_1: P>R$.
\end{enumerate}
The threshold of all three tests is fixed at $5 \%$. 
Throughout this section, we call ``page'' a set of 700 consecutive letters
extracted from a text in our corpus. A \emph{random} page is made of
700 hundred consecutive letters starting at some letter whose rank is
chosen uniformly between $1$ and $N-699$, if there are $N$ letters in
the text. We are interested in the number of pages, needed to raise
the probability of rejecting $\mathcal{H}_0$, i.e. the 
\emph{power} of the tests, above $95 \%$. 
\vskip 2mm\noindent
Assume first a Bernoulli model for the alternance of
vowels and consonants. 
The model states that each letter is a vowel with probability
$p$, or a consonant with probability $1-p$, \emph{independently of the
  others}. The probability distribution of the number of vowels in a
sample of size $n$ is binomial with parameters $n$ and $p$.
The determination of $n$ such that the power of the
test is higher than $95 \%$ is a standard calculation.
Using the data of Table \ref{tab:PO2}, the minimal values of $n$ were
computed. Dividing them by $700$, the following values are found,
for the number of pages necessary to raise the power of the test above
$95 \%$.
\begin{enumerate}
\item[[T1\!\!]] at least $88$ pages,
\item[[T2\!\!]] at least $343$ pages,
\item[[T3\!\!]] at least $52$ pages.
\end{enumerate}
\vskip 2mm\noindent
The Bernoulli model is quite unrealistic: vowels and consonants \emph{do
not} alternate independently in a text. A consequence is that
the actual standard deviation of letter counts on a random page is
smaller than the theoretical standard deviations from the Bernoulli
model. To illustrate this, we have extracted 10000 random pages of
each text in our corpus, counted vowels on each page, and returned the
mean and standard deviations of the 10000 proportions. Table
\ref{tab:PO3} shows the results, for the 20 texts: estimated standard
deviations are about $50 \%$ of what they would be, if
vowel occurrences were independent. The 10000 proportions counted on
each text are represented as boxplots on Figure
\ref{fig:PO}. Alberti's observations are indeed verified \emph{on
  average}. Yet, in spite of the rather small standard deviation, the
vowel counts on one given page may have important variations: counting
one single page of poetry and one of orations at random, might easily 
lead to the conclusion that Alberti was wrong. 
\begin{table}
\begin{center}
\begin{tabular}{|rrrrrrrrrr|}
\hline
\multicolumn{10}{|c|}{Poets}\\\hline
CA&JS&LN&ME&OM&PE&SP&ST&TE&VE\\\hline
0.0122&0.0087&0.0091&0.0108&0.0080&0.0106&0.0090&0.0115&0.0143&0.0124\\\hline\hline
\multicolumn{10}{|c|}{Orators}\\\hline
AM&CG&CP&HS&LP&PP&QD&SC&SO&TO\\\hline
0.0115&0.0100&0.0094&0.0095&0.0089&0.0087&0.0081&0.0100&0.0109&0.0099\\\hline
\end{tabular}
\caption{Estimated standard deviations for the proportions of vowels in 10000
  random pages of each text in our corpus. Under the Bernoulli model,
  these standard deviations should be close to 
$\sqrt{0.44(1-0.44)/700}=0.0188$.}
\label{tab:PO3}
\end{center}
\end{table}
\vskip 2mm\noindent
The following experimental setting was chosen. Selecting $k$ random pages from
each text of our corpus, the proportion of vowels was computed for
each page. The 3 tests T1, T2, T3, were then applied to the 
$10\times k$ proportions among poets and the $10\times
k$ proportions among orators, and the 3 p-values were returned. This was
repeated $10000$ times. For each test, the proportion of the $10000$
p-values below $5 \%$, estimates the power of the test. 
Table \ref{tab:tests} reports the estimated powers for the 3 tests,
with $k=1,\ldots,6$ random pages extracted from each text. The first
values of $k$ such that the estimated power is greater than $95 \%$ are:
\begin{enumerate}
\item[[T1\!\!]] $k=3$, i.e. $60$ pages,
\item[[T2\!\!]] $k=6$, i.e. $120$ pages,
\item[[T3\!\!]] $k=2$, i.e. $40$ pages.
\end{enumerate}
As expected, these values are smaller than those predicted by the
Bernoulli model.
\setlength{\tabcolsep}{10pt}
\begin{table}
\begin{center}
\begin{tabular}{|c|ccc|}
\hline
pages&\multicolumn{3}{|c|}{tests}\\\hline
k&T1&T2&T3\\\hline
1&51.14&33.00&71.16\\
2&84.81&58.24&95.42\\
3&96.23&76.61&99.54\\
4&99.04&86.45&99.95\\
5&99.79&92.73&100\\
6&99.92&96.21&100\\\hline
\end{tabular}
\caption{Estimated powers (in percents) of tests T1, T2, T3 with $k$
  random pages extracted from each of the 20 texts. The number of
  pages in each text needed to raise the power above $95 \%$ is
$k=3$ for T1 ($P>7/16$), $k=6$ for T2 ($R>3/7$), and $k=2$ for T3 ($P>R$).}
\label{tab:tests}
\end{center}
\end{table}
\bibliographystyle{decsi}

%
\end{document}